\documentstyle[11pt]{amsart}
\newtheorem{prop}{Proposition}
\newtheorem{lemma}{Lemma}
\newtheorem{claim}{Claim}
\newcommand{\Y}{{Y}_{3}}
\newcommand{\Z}{\tilde{Z}}
\newcommand{\cL}{{\cal{L}}}
\newcommand{\Jac}{\operatorname{Jac}}
\newcommand{\Pic}{\operatorname{Pic}}
\newcommand{\Tor}{\operatorname{Tors}}
\newcommand{\GH}{\operatorname{H}}
\newcommand{\R}{\right. \\ \left.}
\begin{document}
\title{A surface of general type with
$\rm{p_g}$
\boldmath{$=$} $\rm{q}$
\boldmath{$=0,$}
{\boldmath{${K}^{2}=1$.}
}}
\author{Caryn Werner}
\maketitle

\section{Introduction}
In this paper we construct a
minimal surface \(X\) of general type with
\(\rm{{p}_{g}}=0,\rm{q}=0, {K}^{2}=1,\)
and \( \Tor X \cong \Bbb{Z}/{2}\).
In \cite{Ca}, Campedelli noted that if a degree ten
plane curve could be found having certain singularities,
a double plane construction would yield a surface with
$\rm{{p}_{g}} = \rm{q} =0.$
In \cite{OP}, Oort and Peters construct such a
double plane and compute the torsion group of
their surface to be $\Bbb{Z}/{2}$; however we will
show that the torson group is actually
$\Bbb{Z}/{4}$.
(Weng Lin has also constructed  surfaces with
$\rm{{p}_{g}}=\rm{q}=0,K^{2} =1$ using double covers, but
I do not believe his results are published.)

For minimal surfaces of general type with
$\rm{{p}_{g}}=\rm{q}=0$ and ${K}^{2}=1$
it is known that \( \left| \Tor X \right| \leq 5.\)
(See for example \cite{Do}.)
By writing down generators for the
pluricanonical rings, Reid \cite{R} has described surfaces with  torsion
of order three, four, and five.
Barlow \cite{Ba1,Ba2} has
constructed surfaces with torsion of order two and four,
as well
as a simply connected surface.

The  double plane
constructions  we call {\em Campedelli surfaces},
while {\em numerical Godeaux surfaces} are
minimal surfaces of general type with
the invariants $\rm{p_{g}}=\rm{q}=0, {K}^{2}=1.$
Here \(\rm{{p}_{g}} =
\dim \GH^{0}\left( X,\cal{O}_{X}\left( {K} \right)\right)
= \dim \GH^{2} \left(X, \cal{O}_{X} \right),
\, \rm{q} = \dim \GH^{1} \left( X, \cal{O}_{X} \right), \)
and \( {K}^{2} = {K} \cdot {K} \) is the
self-intersection number of the canonical class \({K}.\)
Write $h^{i}\left(D\right)=
\dim_{\Bbb{C}} H^{i}\left(X,\cal{O}_{X}\left(D\right)\right)$
for $D$ a divisor on $X,$
\( \Tor X\) for the torsion subgroup of the
Picard group of $X,$  \( \equiv\) to represent linear
equivalence of divisors, and \(\left| D \right|\)
for the complete linear system of a divisor class $D.$

\section{The double plane construction.}
Let $D$ be a degree ten plane curve with
an ordinary order four point at $p,$  five
infinitely near triple points at $p_{1}, \dots, p_{5},$
and no other singularities.
An infinitely near triple point
refers to a triple point which remains of order three after
the plane is blown up at this point, so that all three tangent
directions of $D$ coincide.
We assume that each triple point becomes ordinary
after one blow up.
Assume further that
the six singular points do not lie on a conic,
and that the system of plane quartics with
 double point at $p$ and through each $p_{i}$ with
 the same tangent direction as $D$ is exactly
 a pencil.


Let \( \sigma_{1} :{Y}_{1} \rightarrow \Bbb{P}^{2} \)
 be the blowup of \( \Bbb{P}^{2} \) at ${p},$ and let
 ${E}=\sigma_{1}^{-1}({p})$ be the exceptional curve on ${Y}_{1}.$


The total transform of ${D}$ is \( \sigma_{1}^{\ast}
 \left( D \right) = \Bar{{D}} + 4 {E},\)
 where \( \Bar{{D}} \) is the proper transform of $D$.
  Set \( {D}_{1}=\Bar{{D}}
 = \sigma_{1}^{\ast} \left( D \right) -4{E}.\)

 Now let \( \sigma_{2} : {Y}_{2} \rightarrow {Y}_{1} \)
 be the blowup of ${Y}_{1}$ at \( {p}_{1}, \dots, {p}_{5}. \)  With
 \( {E}_{i} =\sigma_{2}^{-1} \left( {p}_{i} \right), \)
 the total transform of \( {D}_{1} \) is
 \( \Bar{{D}_{1}} + 3 \displaystyle{\sum_{1}^{5} {E}_{i},} \)
where $\Bar{{D}_{1}}$ is the proper transform of
$D_{1}.$
Set
 \[ {D}_{2} = \Bar{{D}_{1}} + \sum_{1}^{5} {E}_{i}
 \equiv \sigma_{2}^{\ast}\left( {D}_{1} \right) -2 \sum_{1}^{5} {E}_{i}, \]
which is the reduced divisor consisting of the proper transform of the degree
ten curve, together with the five exceptional curves $E_{i}.$
 As each ${p}_{i}$ is an infinitely near triple point of
 the original branch curve, ${D}_{2}$ has  an order four point on
 each ${E}_{i}.$

 Let \( \sigma_{3} : {Y}_{3} \rightarrow {Y}_{2} \)
 be the blowup of each of these quadruple points, and
 let \( {F}_{1}, \dots, {F}_{5} \) be the corresponding
 exceptional divisors.
We will write ${E}_{i}$ to denote both the
 exceptional curve on $Y_{2}$ and its proper transform
on $Y_{3}$ ( and similarly for ${E}$)
 so that
 \( \sigma_{3}^{\ast} \left( {E}_{i} \right)
 = {E}_{i} + {F}_{i}.
 \)


The total transform of
 ${D}_{2}$ is
 \( \displaystyle{\Bar{{D}_{2}} + 4 \sum_{1}^{5} {F}_{i},} \)
where $\Bar{{D}_{2}}$ is the proper transform of
$D_{2};$ set
 \[
 \begin{array}{cl}
 {B} &= \Bar{{D}_{2}}  \\
 &\equiv \sigma_{3}^{\ast} \left( {D}_{2} \right) -4 \sum {F}_{i} \\
 &= \sigma_{3}^{\ast} \left( \sigma_{2}^{\ast} \left( \sigma_{1}^{\ast}
 \left( {D} \right) -4 {E} \right)
 -2 \sum {E}_{i} \right) -4 \sum {F}_{i} \\
 &= \sigma^{\ast} \left( {D} \right) -4 {E} -2 \sum {E}_{i}
 -6 \sum {F}_{i}
 \end{array}
 \]
where $\sigma = \sigma_{1} \circ \sigma_{2} \circ \sigma_{3}.$
If $H$ represents the pullback of the hyperplane class in
\( \Bbb{P}^{2},\)
 then
\[
 {B}
 \equiv 10H -4{E} -2 \sum_{1}^{5}{E}_{i} -6\sum_{1}^{5}{F}_{i}
 = 2 \cal{L}
 \]
 where
 \[
 \cal{L} \equiv 5H -2E-\sum {E_{i}} -3 \sum{F_{i}};
 \]
\( {B} \) is now a non-singular even curve on the surface
 \( {Y}_{3} . \)
The canonical divisor  on \({Y}_{3}\) is
\[ {K}_{{Y}_{3}} \equiv
\sigma^{\ast}\left(K_{\Bbb{P}^{2}}\right) +{E} +\sum_{1}^{5} {E}_{i}
+2\sum_{1}^{5}{F}_{i} \equiv
-3H +{E} +\sum_{1}^{5} {E}_{i}
+2\sum_{1}^{5} {F}_{i} .\]

 Let \( \pi : {X} \rightarrow {Y}_{3} \) be
 the double cover of ${Y}_{3}$ branched at \( {B}. \)
Then
 \[
 \begin{array}{ll}
 {K}_{{X}} &\equiv \pi^{\ast} \left( {K}_{{Y}_{3}} + \cal{L} \right) \\
 &\equiv \pi^{\ast} \left( 2 {H} - {E} - \sum {F}_{i} \right)
 \end{array}
 \]
 and \( {{K}_{X}}^{2} = 2 {\left[ 2 {H} -{E}
 -\sum {F}_{i} \right]}^{2} = 2 \left( 4 -1 +5\left(-1 \right) \right) = -4. \)
 Since each ${E}_{i}$ is part of
 the branch locus and \( {{E}_{i}}^{2} = -2 \)
 on \( Y_{3}, {\left[ \pi^{-1} \left( {E}_{i} \right) \right] }^{2} = -1.\)
Let
\({X} \rightarrow \tilde{{X}} \) be the map contracting
these five \( \left( -1 \right) \) curves.  Then
\[ \left(
{K}_{\tilde{{X}}} \right)^{2}
=
\left( {K_{X}}\right)^{2} +5
= 1.\]

\begin{prop}
\( \tilde{\bf{X}} \) is  a  minimal surface of general type
with
\( \rm{p_{g}} =0, \rm{q} =0,\mbox{and} \;{K}^{2}
=\chi = 1.\)
\end{prop}

Since  \(\tilde{\bf{X}}\) is obtained from \(\bf{X}\) by
blowing down five exceptional curves, we can compute
\( \rm{{p}_{g}}\) and \(\rm{q}\) for the surface  \(\bf{X}.\)
To compute these invariants, we will use
the following.
\vskip .25 cm

\noindent{\bf Projection formula.}
\begin{sl}
Let \( \pi: X \rightarrow Y\) be a double cover
branched along a smooth curve $B \equiv 2\cal{L}.$
For any divisor ${\cal A}$ on $Y,$
\[
\pi_{\ast}\cal{O}_{X}\left( \pi^{\ast}\cal{A}\right)
\cong
\cal{O}_{Y}\left(\cal{A}\right)
\oplus \cal{O}_{Y}\left( \cal{A}-\cal{L}\right).
\] \end{sl}

In particular,
\[
\pi_{\ast} \cal{O}_{X}\left( nK_{X}\right) \cong
\cal{O}_{Y}\left( nK_{Y}+n\cal{L}\right) \oplus
\cal{O}_{Y}\left( nK_{Y} + \left(n-1\right)\cal{L}\right).
\]
Therefore in our example,
 \[
 \GH^{0}\left(\cal{O}_{{X}}\left( {K}_{X} \right)\right) \cong
 \GH^{0}\left(\cal{O}_{{\Y}} \left( {K}_{{\Y}}+\cal{L} \right)\right)
 \oplus \GH^{0}\left(\cal{O}_{\Y} \left( {K}_{{\Y}} \right)\right),
\]
so that
\[
 p_{g} \left( {X} \right)=
 h^{0}\left( \cal{O}_{X}\left({K_{X}} \right) \right)
 = h^{0} \left( {K}_{{\Y}} + \cal{L} \right) +
 p_{g} \left( {\Y} \right).
 \]
 Since  \( p_{g} \left( {\Y} \right) = p_{g} \left( \Bbb{P}^{2} \right) =0, \)
  \( p_{g} \left( {X} \right) = h^{0} \left(K_{\Y}+\cal{L}\right).\)
  The space $\GH^{0}\left( K_{\Y}+\cal{L}\right)$
  corresponds to the  the linear system
  \(\left| 2 {H} - {E} -\sum {F}_{i} \right| \)
  of conics through \( {p}, {p}_{1}, \dots, {p}_{5}, \)
  so  \( p_{g} \left( {X} \right) = 0.\)
Also
\[
\GH^{0} \left( \cal{O}_{{X}} \left( 2 {K}_{{X}} \right) \right)
= \GH^{0} \left( \cal{O}_{{\Y}} \left( 2 {K}_{{\Y}} + 2 \cal{L} \right)
\right) \oplus \GH^{0} \left( \cal{O}_{{\Y}} \left(
2 {K}_{{\Y}} + \cal{L} \right) \right).
\]
Since
\( \displaystyle{
2K_{\Y} + \cal{L} \equiv -H +\sum_{1}^{5}{\left( E_{i} +F_{i}\right)},
} \)
 \( \GH^{0}\left( 2K_{\Y} +\cal{L}\right) =0.\)
The divisor $\sum E_i$ is a fixed part of
the linear system
\[ \left| 2 {K}_{{\Y}} + 2 \cal{L} \right| =
\left| 4 {H} - 2 {E} -2 \sum_{1}^{5} {F}_{i} \right|
\]
since $E_i \cdot F_i =1$ and $E_i^2 =-2$;
the difference
\( \left| 4H-2E-2\sum F_i -\sum E_i \right| \)
corresponds to quartics in \( \Bbb{P}^{2} \) with a double point
at \( \it{p} ,\) through each \( \it{p}_{i} \) with the same
tangent direction as the branch curve.
By assumption
this system is
a pencil,
thus
\[
\rm{P}_{2} = h^{0}\left( 2 {K}_{{X}} \right) =
\dim \GH^{0} \left( \cal{O}_{{X}} \left( 2 {K}_{{X}}
\right) \right) = 2.
\]

Suppose $S$ is the minimal model of $\tilde{X};$ then
$\rm{P}_{2}\left( S \right) = 2$ and
${K_{S}}^{2} \geq {K_{\tilde{X}}}^{2} =1,$
so $S$ is of general type (see for example \cite{BPV}).
But $2=\rm{P}_{2} = \chi + K_{S}^{2} =1 +K_{S}^{2},$
so $K_{S}^{2} =K_{\tilde{X}}^{2}$ and $S=\tilde{X}.$
Thus \( \tilde{X} \) is minimal of general type with
$K^2 =1$;
since $q \leq p_g $ \cite[\S{3}.1,lemma 3]{Do},
\( p_{g} = q =0  \).

\section{The branch curve {D}}
\newcommand{\pF}{F_{\phi}}
\newcommand{\pset}{p,p_{1},\dots,p_{5}}

To construct a plane curve of degree ten with the
necessary singularities, we will find an octic and
a conic as follows.

We wish to find an octic $C$ with one order four point, one infinitely
near triple point, and four tacnodes, where a tacnode refers to a
double point which remains double after one blowup.
Furthermore we want these tacnodes to lie on a conic $Q$ with the
same tangent direction, so that the octic and conic will
still intersect after the plane is blown up at these points.


Let $F$ be a homogeneous polynomial of degree eight in three variables
defining an octic
${C}$
in $\Bbb{P}^{2}.$
After imposing  an order four point
at ${p} = [1:0:0]$ and an infinitely near triple point
at ${p}_{1} = [0:1:0],$
$F$ has $23$ free coefficients.

Let
\[
\begin{array}{cl}
\gamma : \Bbb{P}^{1} &\rightarrow {Q}  \\
\left[ s : t \right] &\rightarrow \left[ as^2+bst+ct^2 :
ds^2+est+ft^2 : gs^2+hst+it^2 \right]
\end{array}
\]
be a parametrization of a conic \( {Q} \) in \( \Bbb{P}^{2}, \)
where \( a,b,c,d,e,f,g,h,i \) are variables over $\Bbb{C}.$

Set
\[
\begin{array}{cl}
{p}_{2} &= \gamma([0:1])  \\
{p}_{3} &= \gamma([1:0])  \\
{p}_{4} &= \gamma([1:1])  \\
{p}_{5} &= \gamma([-1:1]).
\end{array}
\]

The condition that $F$ have a double point at $p_{i}$ can
be expressed by requiring the three partial derivatives
of $F$  at $p_{i}$ to vanish, thus a double point
is three linear conditions on the coefficients of $F;$
a tacnode at a given point with a designated tangent direction
puts six conditions on $F,$ while
a cusp at a given point with a given tangent direction
is five linear conditions on the
coefficients.
If we impose  tacnodes tangent to $Q$ on the octic at
\( {p}_{2} \) and
\( {p}_{3} \), this gives twelve linear relations
on the coefficients of $F.$
Imposing  cusps tangent to $Q$ at \( {p}_{4} \) and \( {p}_{5} \)
gives ten more relations;
solving these
gives an octic whose coefficients are polynomials in
\( a,b,c,d,e,f,g,h,i .\)
Imposing the  conditions that \( {p}_{4} \)
and \( {p}_{5} \) be tacnodes of \( {C} \)
gives two more linear relations in the coefficients of $F,$
and therefore two higher degree polynomials in
\( a,b,c,d,e,f,g,h,i .\)
In solving these two relations for
\(a,b,c,d,e,f,g,h,\mbox{and} \, i\) we hope to obtain
an irreducible polynomial \(F\) over
\( \Bbb{C}, \) and thus an octic plane curve as desired.

We use Maple to compute the equations for these conditions
on $F,$ and to find the coefficients.
Let \( \left\{ \rm{A}_{j} \right\}_{1}^{22} \) be the equations
corresponding to these  conditions on \( F; \)
the \( \rm{A}_{j} \) are linear in the coefficients of $F.$

Form the matrix {\bf M} generated by the
\( \rm{A}_{j} \) where {\bf M}\(_{i,j} \)
is the coefficient in \( \rm{A}_{j} \)  of the \( ith \) coefficient of
\(F. \)  Then {\bf M} is a $22 \times 23$ matrix, and if we set
{\bf D}\(_{j}\) to be the determinant of the matrix obtained from {\bf M}
be deleting the $j$th column, we have
{\bf M}\(((-1)^{j}\){\bf D}\(_{j}) =0, \)
so that setting the $jth$ coefficient of $F$ to be
$(-1)^{j}${\bf D}$_{j}$ gives the desired octic.

In order for Maple to compute these determinants
quickly enough, we first set
\(a=e=g=i=1\) and
\(b=d=h=0\) in the parametrization of ${Q}.$
This reduces the number of free parameters in this
problem to two, namely $c$ and $f;$  since we will end up imposing two
non-linear conditions on the remaining parameters,
there is still the possibility of a non-degenerate
solution.
After finding
the determinants {\bf D}$_{j}$, the coefficients of $F$ become polynomials
in $c$ and $f.$  Imposing
the final two conditions on the octic, Maple finds several
degenerate solutions, where the octic splits into several
curves of smaller degree, and thus has more singularities, and
also a solution for $c$ and $f$ giving an
octic which we will show has  the desired properties.

The branch curve $D$ is defined by
the equations for the conic and the octic, which
are both polynomials in three variables over
\( \Bbb{Z} \left[ \alpha, \beta, \delta \right] \)
where
\[
\begin{array}{ll}
\alpha &= \sqrt{17} \\
\beta &= \sqrt{21 +5 \sqrt{17}}  \\
\delta &= \sqrt{5 + \sqrt{17}}.
\end{array}
\]

The polynomial  defining the conic $\cal{Q},$ which is
given parametrically by
$\gamma,$ is
\[
\begin{array}{l}
\left(
9\,\alpha \beta
+90\,\alpha +81\, \beta
+234\,  \right) x^{2}
+\left(
+176\,\alpha \beta
+1568\,\alpha +1200\, \beta+5920  \right) y^{2}
\\
+ \left(
57\,\alpha \beta +258\, \alpha +129\, \beta
+1170 \right) z^{2}
+\left(
-48\,\alpha \beta \delta  -168\, \alpha \delta -48\, \beta \delta
-936\,\delta  \right) xy
\\
+\left(
-66\, \alpha \beta
-348\,\alpha
-210\,\beta -1404 \right) xz
+\left(
48\, \alpha \beta \delta
+168\, \alpha \delta
+48\,\beta \delta +936\, \delta
\right) yz
\end{array}
\]

The octic ${C}$ is defined by
$F=0$ where $F$ is
\[
\begin{array}{c}
24\,\left (14408408592\,x^{4}y^{2}z+50076004923\,x^{4}z^{3}+
14182182144\,x^{3}y^{4}
\R
+219953469600\,x^{3}y^{2}z^{2}-363210576777\,x^
{3}z^{4}-1093337332608\,x^{2}y^{2}z^{3}
\R  +831133690121\,x^{2}z^{5}
+
858975454416\,xy^{2}z^{4}-772939669603\,xz^{6}
\R   +254940551336\,z^{7}
\right )y \alpha\beta\delta
\\
+\left (-72389196288\,x^{4}y^{4}-3335393797632\,x^{4}y^{2}
z^{2}-1065820046526\,x^{4}z^{4}
\R
-8342111361024\,x^{3}y^{4}z
+
3945428471808\,x^{3}y^{2}z^{3}
+10184161263912\,x^{3}z^{5}
\R
+
20168534212608\,x^{2}y^{4}z^{2}+53110876008192\,x^{2}y^{2}z^{4}-
32270723397636\,x^{2}z^{6}
\R
-100932292129536\,xy^{2}z^{5}+38252243189640
\,xz^{7}+47211381447168\,y^{2}z^{6}
\R -15099861009390\,z^{8}\right )
\alpha \beta +
\end{array}
\]

\[
\begin{array}{c}
144\,\left (15490159728\,x^{4}y^{2}z+53840161671\,x^{4}z^{3}+
15251365120\,x^{3}y^{4}
\R
+236488232416\,x^{3}y^{2}z^{2}-390512591333\,x^
{3}z^{4}-1175521621376\,x^{2}y^{2}z^{3}
\R +893608694925\,x^{2}z^{5}
+
923543229232\,xy^{2}z^{4}-831040262535\,xz^{6}
\R  +274103997272\,z^{7}
\right )y \alpha \delta
\\
+24\,\left (59398585488\,x^{4}y^{2}z+206468708787\,x^{4}z^{
3}+58496365824\,x^{3}y^{4}
\R
+906899335584\,x^{3}y^{2}z^{2}-1497555836337
\,x^{3}z^{4}-4507944789888\,x^{2}y^{2}z^{3}
\R +3426852351793\,x^{2}z^{5}
+
3541646868816\,xy^{2}z^{4}-3186912029723\,xz^{6}
\R +1051146805480\,z^{7}
\right )y \beta \delta
\\
+\left (-466877917440\,x^{4}y^{4}-21516582641184\,x^{4}y^{2
}z^{2}-6875617032333\,x^{4}z^{4}
 \R  -53815549731840\,x^{3}y^{4}z
+25451977456512\,x^{3}y^{2}z^{3}+65698123816692\,x^{3}z^{5}
\R  +
130107481479168\,x^{2}y^{4}z^{2}
+342618718898784\,x^{2}y^{2}z^{4}-
208178748305934\,x^{2}z^{6}
\R -651115201689280\,xy^{2}z^{5}  +
246765593291124\,xz^{7}+304561087975168\,y^{2}z^{6}
\R  -97409351769549\,z^
{8}\right ) \alpha
\\
+\left (-298344909312\,x^{4}y^{4}-13752106145280\,x^{4}y^{
2}z^{2}-4394491299054\,x^{4}z^{4}
\R -34395989170176\,x^{3}y^{4}z
+
16267404201984\,x^{3}y^{2}z^{3}
+41990375439720\,x^{3}z^{5}
\R +
83157067186176\,x^{2}y^{4}z^{2}+218981688444672\,x^{2}y^{2}z^{4}-
133055599121316\,x^{2}z^{6}
\R
-416154499902208\,xy^{2}z^{5}+
157718037119688\,xz^{7}+194657513400832\,y^{2}z^{6}
\R  -62258322139038\,z^
{8}\right ) \beta
\\
+48\,\left (191605550544\,x^{4}y^{2}z+665965983645\,x^{4}z
^{3}+188641373952\,x^{3}y^{4}
\R   +2925195141792\,x^{3}y^{2}z^{2}
-4830373865031\,x^{3}z^{4}-14540399432832\,x^{2}y^{2}z^{3}
\R  +
11053328983807\,x^{2}z^{5}
+11423598740496\,xy^{2}z^{4}-10279400307101
\,xz^{6}
\R  +3390479204680\,z^{7}\right )y \delta
\\
-1925078503680\,x^{4}y^{4}-
88715185482528\,x^{4}y^{2}z^{2}-28348893570645\,x^{4}z^{4}
\\
-
221886790124544\,x^{3}y^{4}z+104941198124928\,x^{3}y^{2}z^{3}+
270880301999124\,x^{3}z^{5}
\\
+536446841591808\,x^{2}y^{4}z^{2}+
1412653204386144\,x^{2}y^{2}z^{4}-858342969385662\,x^{2}z^{6}
\\  -
2684616751584448\,xy^{2}z^{5}
+1017440607056532\,xz^{7}+
1255737534555904\,y^{2}z^{6}
\\-401629046099349\,z^{8}.
\end{array}
\]

The resulting singular points of the branch curve are
\[
\begin{array}{ll}
p &= [1:0:0] \\
p_{1} &= [0:1:0] \\
p_{2} &=
[{{10
+{4\,\alpha}+{4\,\beta}}}:{ 3{\delta}}:6]
\\
p_{3} &=
[1:0:1] \\
p_{4} &=
[{{ {16}
+{4\,\alpha}+{4\,\beta}}}:{3{\delta}}+6:12]
\\
p_{5} &=
[{{16}
+{4\,\alpha}+{4\,\beta}}:{ 3{\delta}}-6:12].
\end{array}
\]

We need to check that the branch curve $D$ has no singularities outside the set
$\left\{ \pset \right\}.$
Since $F$ is  a polynomial over
the complex numbers, Maple is unable to
quickly check that the octic has no other singularities,
so we use Macaulay
to check the smoothness of ${C}$ outside of the set
\( \left\{p,p_{1},\dots,p_{5} \right\}.\)
As Macaulay only makes computations over finite
fields,  we first find a prime number $P$ where
$\alpha,\beta,\mbox{and} \, \delta$
exist mod $P,$ so that
we can map $F$ to a polynomial over $\Bbb{Z}/{P}.$

To check that ${C}$ has no singularities other than
at the points \( {p}, {p}_{1}, \dots, {p}_{5}, \)
consider the map
\(
\phi: \Bbb{Z}[\alpha,\beta,\delta] \rightarrow \Bbb{Z}/30047
\)
given by sending
\[
\begin{array}{ll}
\alpha &\rightarrow 20452 \\
\beta &\rightarrow  6941 \\
\delta &\rightarrow 27962;
\end{array}
\]
mapping $F$ to $F_{\phi}$
we obtain
\[
\begin{array}{c}
24082\,x^{4}y^{4}+3438\,x^{4}y^{3}z+4775\,x^{4}y^{2}z^{2}+29499\,x^{4}
yz^{3}+12698\,x^{4}z^{4}
\\
+29927\,x^{3}y^{5}+14121\,x^{3}y^{4}z+17243\,x
^{3}y^{3}z^{2}+3139\,x^{3}y^{2}z^{3}+8704\,x^{3}yz^{4}+80\,x^{3}z^{5}
\\ +
28712\,x^{2}y^{4}z^{2}+10654\,x^{2}y^{3}z^{3}+12817\,x^{2}y^{2}z^{4}+
8239\,x^{2}yz^{5}+5515\,x^{2}z^{6}
\\
+28759\,xy^{3}z^{4}+7372\,xy^{2}z^{5
}+19696\,xyz^{6}+28079\,xz^{7}
\\
+1944\,y^{2}z^{6}+24003\,yz^{7}+13722\,z^{8}.
\end{array}
\]

\begin{claim}
$\pF$ has no singularities other than at $\pset.$
\end{claim}

First we have Macaulay compute $\Jac \pF,$ the Jacobian
ideal of the octic  generated by
\(\displaystyle{ {{\partial{\pF}} \over {\partial{x}}},
 {{\partial{\pF}} \over {\partial{y}}}, \mbox{and}
 {{\partial{\pF}} \over {\partial{z}}}, }\)
 and the ideal $\cal{I}$ associated to the points
\( {p},{p}_{1},{p}_{2},p_{3},p_{4},p_{5}.\)
Since the zeros of
\( \Jac\left( \pF \right) \) are precisely the
singular points of the octic, the zeros of
the saturation of \( \Jac \left( \pF \right) \)
by \( \cal{I} \) are any singularites other than
at the zeros of $\cal{I}.$
Macaulay computes
\[
\left( \Jac\left( {\pF} \right) : \cal{I}^{\infty} \right) =
\left\{ g : g \cal{I}^{n} \subset \Jac \pF \,
\mbox{for some n} \right\}
 = \left( 1 \right),
\]
thus there are no zeros of
$\Jac \pF$ other than at the points
\(p, p_{1},\dots,p_{5}\) and therefore no
other singularities of $C_{\phi}.$

\begin{claim}
To check that $F$ has no singularities outside the
set $\left\{ p,p_{1}, \dots, p_{5} \right\},$ it suffices
to check this for the polynomial $F_{\phi}$
over \( \Bbb{Z}/{30047}.\)
\end{claim}

It is easy to check, using Maple,
that $C_{\phi}$ has an ordinary quadruple point at $p$,
and after one blow up, the triple point at $p_{1}$ and the
double points at $p_{2},\dots,p_{5}$ become ordinary.
Since $C$ maps to $C_{\phi}$, the same is true for the
singularities on $C$.

Maple is not reliable about completely factoring polynomials
in many variables; hence Maple cannot check directly
that $C$ is irreducible.
Hence we fall back on a more case-by-case analysis.

Maple {\em can} check that a given polynomial divides another;
and so one can use Maple to conclude that
$Q$ is not a component of $C$.

Next note that since $\deg Q = 2$ and $\deg C = 8$, $Q \cdot C =16$.
We know that $C$ and $Q$ meet four times at each $p_{i},
i=2,\dots,5$; thus $Q$ cannot meet any component of $C$
at any other point.

We check that none of the tangent lines to $C$ at any $p_i$
are contained in $C$.  If any other line was a component
of $C$, say $C=\ell G$, then $\left( G \cdot Q \right) =14$;
however $G$ must meet $Q$ four times each at $p_2,\dots,p_5$,
so no line can be contained in $C$.

Suppose a conic $G$ is a component of $C$.  Then $G$ must
meet $Q$ at two of the four points, say $p_{i}$ and $p_{j}$,
with the proper tangent directions, to multiplicity two.

{Case 1.}  If $C$ breaks up into $G$ and an irreducible sextic $S$,
then $S$ must have at least a triple point at $p$ and
tacnodes at $p_{k},p_{l}$ for $l,k \neq i,j$; since there
is no conic through $p_{1},p_{i},p_{j}$ with the required
tangent directions, $S$ must have an infinitely near triple
point at $p_{1}$.  But these conditions would drop
the genus of $S$ by $13$, while an irreducible degree six curve has
genus $10$, so no such sextic exists.

{Case 2.} If $C$ breaks up into two conics $G$ and $H$
and a degree four part,
then $G$ meets $Q$ at $p_{i},p_{j}$,
$H$ meets $Q$ at $p_{k},p_{l}$,
so neither conic can pass through $p_{1}$.
Therefore the degree four part of $C$ would have
to have an infinitely near triple point, which is impossible
(even for a reducible quartic).

{Case 3.}  If $C$ is composed of a conic $G$ and two cubics $S_{1},S_{2}$,
then one of the cubics must have a
tacnode at $p_{1}$, which is impossible.

Thus the octic $C$
cannot contain either a line or a conic as a component.
We can conclude therefore that if $C$ does split,
it splits into at most two components
(of degrees $3$ and $5$ or $4$ and $4$).

Suppose $C$ is composed of a cubic $G$ and a quintic $S$,
both of which are irreducible.
The arithmetic genus of $S$ is six, and $S$ must have
at least a double point at $p$ and a tacnode at $p_{1}$,
which together drop the genus by three.  Since $Q \cdot G =6$,
$G$ must meet $Q$ at three of the $p_{i}$, thus $S$ must
have a tacnode along $Q$ (at the fourth point)
which drops the genus by two more.
Thus $S$ can have exactly a double point at $p$ and a
tacnode at $p_{1}$, and $G$ must have a double point
at $p$ and pass through $p_{1}$ and three of the
$p_{i}$ with the necessary tangent directions.
But no such cubics exist, as can be checked using Maple;
(this gives $11$ linear conditions on the cubic, and Maple checks that this
linear system has no solutions).

Next, suppose $C$ is composed of two irreducible quartics $G$ and $S$.
Then one of the quartics, say $G$, must have a
tacnode at $p_{1}$ and pass through $p$.  Also $G$ must
meet $Q$ along $p_{2},\dots,p_{5}$.  But these are all linear conditions
on the quartic, and again Maple can be used to check that there are
no such quartics.

Thus the octic $C$ is irreducible.

Since $C$ is irreducible, we can compute the arithmetic genus
to be
$\left( \begin{array}{c} 7 \\ 2 \end{array} \right) = 21;$
after blowing up a point of multiplicity $n,$ the
genus of the proper transform goes down by
$\left( \begin{array}{c} n \\ 2 \end{array} \right).$
After resolving the singularities of $C$ at $\pset,$
the resulting curve has genus equal to
\[
21 - \left( \begin{array}{c} 4 \\ 2 \end{array} \right)
- 2 \left( \begin{array}{c} 3 \\2 \end{array} \right)
-8 \left( \begin{array}{c} 2 \\ 2 \end{array} \right) =
1, \]
thus $C$ can have at most one more singularity of multiplicity two.

We will now prove that $C$ has no other singularities
than the known ones at $p$ and $p_1,\ldots,p_5$.
The curve $C$ is defined over the field
$K={\bf Q}(\alpha,\beta,\delta)$,
as is its  strict transform $\bar C$
after resolving the singularities at $p,p_1,\dots,p_5$.
Suppose that $\bar C$ is singular;
since it can have at most one singularity,
the coordinates of this singular point are then invariant
by the action of the Galois group of the algebraic closure of $K$ over $K$,
hence lie in $K$.
Thus the normalization $\tilde C$ of the curve $\bar C$ is defined over $K$.
Since the genus of $\tilde C$ is $0$,
its anti-canonical map induces an isomorphism, defined over $K$,
onto a smooth conic in ${\bf P}^2_K$.
Since the curve $\tilde C$ has a rational point over $K$
(namely the eighth point of intersection of the line $z=0$ with the curve $C$:
this line meets $C$ four times at $p$,three times at $p_1$,
and then once at a point with coordinates in $K$),
the projection from this point yields an isomorphism defined over $K$
between $\tilde C$ and ${\bf P}^1_K$.
By composing with the map $\tilde C\rightarrow C$ (also defined over $K$),
we obtain a parametrization $\psi:{\bf P}^1_K\rightarrow C$ defined over $K$;
by clearing denominators we can
take $\psi$ to be defined over
${\Bbb Z} \left[ \alpha, \beta, \delta \right]$.
Since ${\Bbb Z} \left[ \alpha,\beta, \delta \right]$ maps to
${\Bbb Z}/{30047}$,
we get a map ${\Bbb P}^{1}_{{\Bbb Z}/{30047}} \rightarrow C_{\phi}$.
Thus $C_{\phi}$ is rational over ${\Bbb Z}/{30047}$,
so the genus is zero and the genus of $C_{\phi}$
over the algebraic closure of ${\Bbb Z}/{30047}$ is also zero.

But Macaulay can and does
check that $C_{\phi}$ has no other singularities
in the algebraic closure of the finite field;
so $\bar{C_{\phi}}$ is smooth and its genus must be one
(using the genus formula, which is essentially adjunction).
Therefore the genus of $\bar{C}$ must be one as well,
which gives a contradiction.
Hence $C$ can have no other singularities.

We can also use Maple to check that the system of
quartics with a double point at $p,$ through each
$p_i$ with the necessary tangent direction
is a pencil;
thus $C$ and $Q$ give a degree ten curve as needed.

\section{The torsion group of \( \tilde{X} \)}

The following lemma will show that the torsion group is
non-trivial.
\begin{lemma}(Beauville \cite{B})
Let  \( Y\) be a smooth surface  with \(\Tor(\Pic(Y)) =0,\)
\( \left\{ C_{i} \right\}_{i \in I} \) a collection of
smooth disjoint curves on \( Y ,\) and
\( \pi : X \rightarrow Y \) a connected  double cover branched
along  \( \cup_{i \in I}C_{i}.\)  Define a map
\[ \varphi : {\Bbb{Z}/{2}}^{I} \rightarrow \Pic Y
\otimes \Bbb{Z}/{2} \]
by sending \( \sum n_{i} C_{i} \, \) to its class in
\( \Pic Y. \) If \( e =
 \sum_{i \in I} C_{i}  \),
then the group $\Pic_{2} X$ of $2-$ torsion elements in
\( \Pic X\)
is isomorphic to
\(\rm{ker}\left( \varphi \right) / \left(\Bbb{Z}/{2}\right)e.\)
\end{lemma}

If \( \sum_{i \in J} C_{i} \equiv 2A\) for some divisor
\(A, \) where  $J$ is a subset of $I,$ then
the map from \( \rm{ker} \left( \varphi \right) \) to
the $2-$torsion elements
in \( \Pic X \) sends \( \sum_{i \in J} C_{i} \) to
\( \sum_{i \in J} \pi^{-1} \left( C_{i} \right)
- \pi^{\ast} \left(A \right) ; \)
 for components $C_{i}$ of the branch locus
\[ 2 \pi^{-1} \left( C_{i} \right) \equiv \pi^{\ast}
\left( C_{i} \right) , \]
so that
\( \sum_{i \in J} \pi^{-1} \left( C_{i} \right) -
\pi^{\ast} \left( A \right) \) is in \( \Pic_{2} \left( X \right) .\)

Let $\bar{Q}$ be the strict transform of $Q$ on $Y_3$.  Since
$\bar{Q} +\sum E_{i}$ is a sum of components of the branch
locus and \(\bar{Q} +\sum E_{i} \equiv 2 \left( H-\sum F_{i} \right)
\), the lemma shows that the divisor
\[
\pi^{-1} \left( \bar{Q} +\sum_{2}^{5} {E}_{i} \right)
- \pi^{\ast} \left(
H - \sum_{2}^{5} {F}_{i} \right)
\]
has order two in \( \Pic \left( {X} \right) .\)
Thus
\( \Tor \left({X}\right) \) is non-trivial.

For numerical Godeaux surfaces, the torsion group has order
less than or equal to five, and it is known that
$\Bbb{Z}/{2} \oplus \Bbb{Z}/{2}$ does not occur. (See \cite{Do}.)
To determine whether $\Tor X$ is
\( \Bbb{Z}/{2} \) or \( \Bbb{Z}/{4} ,\) we use a base
point lemma due to Miyaoka \cite{M}:
for a minimal Godeaux surface, the number of base points
of $\left| 3 {K} \right| $ is equal to
\[
\# \left\{ T \in \Pic X \, : \, T \neq -T\right\}/2.
\]
Thus if $\, \left| 3 {K} \right|$ has no base points,
the torsion group is  \( \Bbb{Z}/2 \).

Write $\epsilon : X \rightarrow \tilde{X}$ for the
map contracting the $\pi^{-1}\left( E_i \right)$.
Then $3K_X \equiv \epsilon^{\ast}\left( 3K_{\tilde{X}} \right)
+3\sum \pi^{-1}\left(E_i \right)$.
To compute $\left| 3 {K}_{X} \right|,$
first consider the system $\left| 3K_{Y_3}+3\cal{L} \right|$.
The divisor $2\sum E_i$ is fixed in this system; the difference
$\left|
6H-3E-2\sum E_i -3 \sum F_i \right|$
is  the pencil of
sextics with a triple point at $p$ and double points at each
$p_i$ with one tangent direction coinciding with the branch curve.
Set $M=
\pi^{\ast}\left( 6H-3E-2\sum E_i -3 \sum F_i \right)$;
we have $\epsilon^{\ast}\left( 3K_{\tilde{X}} \right) \equiv M+
\sum \pi^{-1} \left( E_i \right)$, so any base point must
either lie on
$\sum \pi^{-1} \left( E_i \right)$ or be a base point of
$\left| M \right|$.

We use Maple to find two  sextics in $M$ and their two points of
intersection.
These two points do not lie on $Q$ or $C$, so there is no base
point of $\left| M \right|$ on the branch curve.

Since $3K_X \equiv \pi^{-1}\left( B \right) + \pi^{\ast}
\left( H-E \right) + 2 \sum \pi^{-1} \left(E_i \right)$,
we also have $\epsilon^{\ast}\left( 3K_{\tilde{X}} \right)
\equiv \pi^{-1} \left( B \right) +\pi^{\ast}\left( H-E \right)
-\sum \pi^{-1}\left( E_i \right) $, so any base point must lie on
the branch curve, away from the divisor
$\sum \pi^{-1} \left( E_i \right)$.

Therefore there are no base points of the tricanonical system.
{}From the Miyaoka lemma, this shows that
$\Tor X \cong \Bbb{Z}/{2}.$

\section{The Oort and Peters example}
In \cite{OP}, Oort and Peters construct a branch curve $B$ from
two conics ${Q}_{1},{Q}_{2}$ and two cubics
${C}_{1},{C}_{2}$ where
\[
\begin{array}{ll}
{Q}_{1} &= y^2+2x^2-2xy-5xz+2yz+3z^{2} \\
{Q}_{2} &= y^2+2x^2+2xy-5xz-2yz+3z^2 \\
{C}_{1} &= y^{2}z+x^3-4x^{2}z+3xz^2 \\
{C}_{2} &= 2y^{2}z-xy^2+4x^{2}z-12xz^{2}+9z^3.
\end{array}
\]
We have
\[
\begin{array}{ll}
( {Q}_{1} \cdot { Q}_{2}) &= {P} + 3 { P}_{1} \\
( {Q}_{1} \cdot { {C}}_{1}) &= 2( {{P}}_{1} + {P}_{2} + {P}_{3} ) \\
({ {Q}}_{1} \cdot { {C}}_{2}) &= 2({{P}} + {P}_{2} + {P}_{3} ) \\
({ {Q}}_{2} \cdot { {C}}_{1}) &= 2( {{P}}_{1} + {P}_{4} + {P}_{5} ) \\
({ {Q}}_{2} \cdot { {C}}_{2}) &= 2({{P}}_{4} + {P}_{5} +{P}) \\
( { {C}}_{1} \cdot { C}_{2}) &= 2({ {P}}_{2} +{P}_{3} +{P}_{4} +{P}_{5} ) +
\infty
\end{array}
\]
where
\[
\begin{array}{ll}
{P} &= [{{3} \over {2}}:0:1] \\
{P}_{1} &= [1:0:1] \\
{P}_{2} &= [ {{3+i \sqrt{3}}\over{2}} :  {{3+i \sqrt{3}}\over{2}}  : 1] \\
{P}_{3} &= [ {{3-i \sqrt{3}}\over{2}}  : {{3-i \sqrt{3}}\over{2}}  : 1] \\
{P}_{4} &= [ {{3+i \sqrt{3}}\over{2}}  : {{-3-i \sqrt{3}}\over{2}}   : 1] \\
{P}_{5} &= [ {{3-i \sqrt{3}}\over{2}}  : {{-3+i \sqrt{3}}\over{2}}   : 1]\\
\infty &= [0:1:0].
\end{array}
\]
In this case the branch curve has two extra
ordinary double points, one at
$\infty,$ and the other which occurs on the second blowup above
${P}_{1},$ since ${Q}_{1}$ and ${Q}_{2}$  intersect
with multiplicity three
at this point.  However these double points do not affect the
invariants of the double plane $Z$ constructed.

Write
\( \pi : Z \rightarrow Y\) for the double cover, where $Y$ is
the blowup of the plane resolving the singularities
of the branch curve;  although the double points
do not affect the computations, we will blow them up to obtain
a smooth branch divisor ${B}$ with
\[
{B} =2 \cal{L} \equiv 10H -4 {E}
-2\sum_{1}^{5} {E}_{i} -6 \sum_{1}^{5} {F}_{i}
-8 {G}_{1} -2 {E}_{6}
\]
where we use the notation for the exceptional curves as above,
with ${G}_{1}$ being the divisor lying above the extra double
point on ${F}_{1}$ and ${E}_{6}$ the exceptional divisor
above $\infty.$
Let $\Z$ be the minimal surface obtained from $Z$ by
blowing down the  $E_i$.

Note that we have the following equivalences of divisors:
\[
\begin{array}{cl}
\bar{Q_1} &\equiv 2H-E-E_1-E_2-E_3-2F_1-2F_2-2F_3-3G_1 \\
\bar{Q_2} &\equiv 2H-E-E_1-E_4-E_5-2F_1-2F_4-2F_5-3G_1 \\
\bar{C_1} &\equiv 3H-\sum_{1}^{5} E_i -2\sum_{1}^{5} F_i -2G_1 -E_6\\
\bar{C_2} &\equiv 3H -2E- \sum_{2}^{5} E_i -2 \sum_{2}^{5} F_i -E_6 \\
K_Y &\equiv -3H +E+\sum_{1}^{5} E_i+2\sum_{1}^{5} F_i +3G_1 +E_6 \\
\cal{L} &\equiv 5H -2E-\sum E_i-3\sum F_i -4 G_1 -E_6.
\end{array}
\]
Let $\cL_1 =2H-E -E_1-2F_1-\sum_{2}^{5} F_i-3G_1$
and $B_1 =2 \cL_1 \equiv Q_1 +Q_2 +\sum_{2}^{5} E_i$;
set $\cL_2 = \cL -\cL_1$ and $B_2 =2\cL_2 \equiv C_1+C_2+E_1$.

It follows from Beauville's lemma that
\[ T =\pi^{-1}\left(B_1 \right) -\pi^{\ast} \left( \cL_1 \right)
\equiv -\pi^{-1}\left( B_2 \right) +\pi^{\ast} \left( \cL_2 \right) \]
is of order two, thus
$\Tor \Z$ is either $\Bbb{Z}/2$ or $\Bbb{Z}/4$.

We will show that $\Tor \Z \cong \Bbb{Z}/{4}.$  Note that it
was previously believed that $\Tor \Z \cong \Bbb{Z}/{2}$
(\cite{Do,OP}).

Oort and Peters use the base point lemma of Miyaoka
to argue that $\Tor \Z$ is $\Bbb{Z}/{2}$; however they
miss a base point of the system $\left| 3K_{\tilde{Z}} \right|$
in their computation.
Again if $\epsilon : Z \rightarrow \Z$ is the map from  $Z$ to
its minimal model, we have
$$\epsilon^{\ast} \left(3K_{\Z} \right) \equiv M
+ \sum \pi^{-1} \left( E_i \right) \equiv
\pi^{-1} \left( B \right) +\pi^{\ast} \left( H-E \right)
-\sum \pi^{-1} \left( E_i \right),$$ where
$M= \pi^{\ast} \left(3K_Y +3\cL -2\sum E_i\right) $. Thus any base point
must lie on $\pi^{-1} \left( B \right)
-\sum \pi^{-1} \left( E_i \right)$ and be a base point of
$\left| M \right|$.
The divisors $\bar{Q_1} +\bar{Q_2} +\bar{Q} +F_1$ and
$\bar{\ell} +\bar{C_2} +\bar{\tilde{Q}}$ are in $\left| M \right|$,
where
$\bar{Q}
\equiv 2H-E-\sum_{2}^{5} \left(E_i +F_i \right)$
is the proper transform of  the conic
\(Q = 2x^{2}-9xz+y^{2}+9z^{2}\)
through $P,P_2,\dots,P_5,$
$\bar{\ell} \equiv H-E-E_1-F_1-G_1$ is the proper transform
of the line $y=0$,
and $\bar{\tilde{Q}} \equiv 2H -\sum_{1}^{5} E_i -2F_1
-\sum_{2}^{5} F_i -2G_1$ is the proper transform of the
conic
\(\tilde{Q} = 3xz-3z^{2}-y^{2} \)
through each $P_i$ where the tangent direction at $P_1$
coincides with that of the branch curve.

The point $[3:0:1]$ lies on the curves
${Q}, \ell, \mbox{and} \, {C}_{1},$
and therefore is a base point of
$\left| \epsilon^{\ast} \left(3K_Z \right) \right|$.
It follows from Miyaoka's result that
$\Tor \Z$ is  $\Bbb{Z}/{4}.$

In \cite{Do}, Dolgachev assumes that there exists an order four divisor
on $Z$ and  gets a contradiction after  finding a fixed part of
the pencil
$\left| 2 {K}_{Z} \right|.$
However his computation of generators for
$\left| 2 {K}_{Z} \right|$ is incorrect.
We find divisors in the system
\[
\begin{array}{ll}
\left| \epsilon^{\ast}\left( 2K_{\Z} \right) \right|
&= \left| 2K_Z-2\sum \pi^{-1}\left(E_i \right) \right| \\
&=
\left| \pi^{\ast}\left(2K_Y +2\cal{L} -\sum E_i \right) \right| \\
&=
\left| \pi^{\ast}\left(4H-2E-\sum E_i -2\sum F_i -2G_1\right) \right|.
\end{array}
\]

This pencil  has generators
\(
y_0=\bar{{Q}_{1}} +\bar{{Q}_{2}} +2{F}_{1}+4{G}_{1}+E_1
\)
and
\(
y_1=\bar{{C}_{2}} + \bar{\tilde{\ell}} +2E_6,
\)
where $\tilde{\ell}$ is the line tangent to the branch
curve at $P_1$;
thus there is no fixed part to this system.

We can also check that $\Z$ has order four torsion by
calculating the bicanonical system
of a double cover
of $\Z.$
Form the double cover  $S$ of $\Z$ branched over $2{T} \equiv 0,$
\( \rho : S \rightarrow \Z. \)
Since there is no ramification, $\rho$ is \'etale over $\Z$.
Also $K_{S} \equiv \rho^{\ast} \left( K_{Z}+T\right)$ and
${K_{S}}^{2} = 2.$
We have already found two sections $y_0$ and $y_1$ in
$\GH^{0} \left( 2K_{\Z} \right)$, and hence two sections
$\rho^{\ast}\left( y_0 \right)$ and
$\rho^{\ast}\left( y_1 \right)$ in
$\GH^{0} \left( 2K_S \right)$.
Since $\rho^{\ast}\left( T \right) \equiv 0$,
we also have
\[
2K_S \equiv \rho^{\ast}\bigl(\pi^{-1}\left(B_1 \right)
+\pi^{\ast}(2H-E-\sum_{2}^{5} \left( E_i +F_i \right)
+G_1 ) \bigr)
\]
and
\[
2K_S \equiv \rho^{\ast}\left(\pi^{-1}\left(B_2 \right)
+\pi^{\ast}\left( H-E-E_1 -F_1-G_1+E_6 \right) \right).
\]
We have seen that the  proper transform
$\bar{Q}$ of the conic $Q$ is in the linear system
\(
\left| \pi^{\ast}\left(2H-E-\sum_{2}^{5} \left( E_i + F_i \right)  \right)
\right| \)
and the proper transform of the line $\ell$
is in
\(
\left| \pi^{\ast}\left( H-E -F_1-G_1 \right)
\right|
\);
set $y_2 = \pi^{-1}\left( B_1 \right) +\pi^{\ast}
\left( \bar{Q}  +G_1 \right)$
and $y_3 = \pi^{-1}\left( B_2 \right) +
\pi^{\ast}\left( \bar{\ell} + E_6 \right)$.
We have
\(
\left(y_{0}-2y_{1} \right)^{2} -y_{2}^{2} +4y_{3}^{2} = 0
\).
This gives a quadratic relation among the four elements of
$\GH^{0}(2{K}_{{S}});$ in fact we obtain a quadric cone
as the bicanonical image of ${S}$.  By \cite{CD},  if the bicanonical
image is a cone then  $\Tor S \cong \Bbb{Z}/{2}$ and
$\pi_{1} \left( S \right) \cong \Bbb{Z}/{2}.$
Since $S$ is a covering space of $\Z$ of degree two,
$\left[ \pi_{1}\left(\Z\right) : \pi_{1}\left(S\right)\right] =2.$
Thus
$\pi_{1} \left( \Z \right)$ is abelian of order four and
$\pi_{1} \left( \Z \right) \cong \Tor \left(\Z \right)
\cong \Bbb{Z}/{4}$.

Department of Mathematics,
Colorado State University,
Fort Collins, CO 80523

E-mail address: werner@@lagrange.math.colostate.edu
\end{document}